\documentclass[11pt]{amsart}
\usepackage{amsmath,amssymb,amsbsy,amsfonts,latexsym,amsopn,amstext,amsxtra,epic, euscript,amscd,indentfirst, tikz}
\usepackage{verbatim} 

\usepackage{hyperref} 
\usepackage{xcolor}
\usepackage{graphicx}
\usepackage[justification=centering]{caption} 

\DeclareFontFamily{U}{mathc}{}
\DeclareFontShape{U}{mathc}{m}{it}%
{<->s*[1.03] mathc10}{}

\DeclareMathAlphabet{\mathscr}{U}{mathc}{m}{it}

\begin{document}

\newtheorem{theorem}{Theorem}
\newtheorem{conjecture}[theorem]{Conjecture}
\newtheorem{proposition}[theorem]{Proposition}
\newtheorem{question}{Question}
\newtheorem{problem}{Problem}
\newtheorem{lemma}[theorem]{Lemma}
\newtheorem{definition}{Definition}
\newtheorem*{definition*}{Definition}
\newtheorem{cor}[theorem]{Corollary}
\newtheorem*{cor*}{Corollary}
\newtheorem*{result*}{Result}
\newtheorem{obs}[theorem]{Observation}
\newtheorem{proc}[theorem]{Procedure}
\newcommand{\comments}[1]{}

\def\Z{\mathbb Z}
\def\Za{\mathbb Z^\ast}
\def\Fq{{\mathbb F}_q}
\def\R{\mathbb R}
\def\N{\mathbb N}
\def\i{\sqrt{-1}}
\def\k{\kappa}

\title[Optimal transport and information geometry]{When Optimal Transport Meets \\ Information Geometry }
 \author[Iowa State University]{Gabriel Khan} 
 \author[University of Michigan]{Jun Zhang}
 \email{gkhan@iastate.edu}
 \email{junz@umich.edu}

\maketitle

\begin{abstract}
Information geometry and optimal transport are two distinct geometric frameworks for modeling families of probability measures. During the recent years, there has been a surge of research endeavors that cut across these two areas and explore their links and interactions. This paper is intended to provide an (incomplete) survey of these works, including entropy-regularized transport, divergence functions arising from $c$-duality, density manifolds and transport information geometry, the para-K\"ahler and K\"ahler geometries underlying optimal transport and the regularity theory for its solutions.
Some outstanding questions that would be of interest to audience of both these two disciplines are posed. Our piece also serves as an introduction to the Special Issue on Optimal Transport of the journal {\it Information Geometry}.    
\end{abstract}

\tableofcontents

\section{Optimal Transport: A Brief Overview}

\subsection{The Monge and Kantorovich formulations}

Optimal transport is a classic area of mathematics which combines ideas from geometry, analysis, measure theory, and probability. Today, it is a thriving area of research, both in the pure and applied settings. Furthermore, it has many practical applications ranging from logistics, economics, computer vision, imaging processing, and many other fields.

Optimal transport was first studied by Gaspard Monge in 1781 \cite{monge1781memoire}. In this work, Monge was interested in finding the most efficient way to construct military fortifications (or other architecture) from a distribution of raw material. To study this topic from a mathematical standpoint, he posed the following question, which is now known as the \emph{Monge problem} of optimal transport.
\begin{problem}[Monge's problem]
A worker is tasked with moving a pile of rubble (\textit{deblais}) into a prescribed configuration (\textit{remblais}). The worker must complete the job, but would like to minimize the effort required, and so move the rubble in the most cost-efficient way possible. How should this be done?
\end{problem}

Although the Monge problem is a natural question and is relatively simple to state, it turns out to be quite difficult to solve, and a satisfactory answer for this problem took over 200 years to develop (see Section \ref{Existence of Monge problem} for details). As a result, much of the work in optimal transport instead uses a framework developed by Kantorovich \cite{kantorovich2006translocation} in the 1940's. To contrast this with the original problem of Monge, we will call this the ``Kantorovich problem" (although it is sometimes known as the  {\em Monge-Kantorovich problem}.)

\begin{problem} [Kantorovich problem of optimal transport]
Let $(X, \mu)$ and $(Y, \nu)$ be two Polish probability spaces and $c:X \times Y \to \mathbb{R} \cup \{\infty \}$ be a lower semicontinuous cost function. The Kantorovich problem seeks to find a coupling $\pi$ of $\mu$ and $\nu$ (i.e. a probability measure on $X \times Y$ whose marginals are $\mu$ and $\nu$, respectively) which achieves the smallest total cost
\begin{equation} \label{Kantorovichproblem}
     \inf_{\pi \in \Pi(\mu,\nu)} \int_{X \times Y} c(x,y) \, d \pi(x,y). 
\end{equation}
\end{problem}

Here, $\Pi(\mu, \nu)$ is the set of all couplings between $\mu$ and $\nu$. Intuitively, $\mu$ is the shape of the original configuration and $\nu$ is the target configuration, and transport plan at a point $x \in X$ is given by the disintegration of the coupling $\pi$. The cost function $c$ gives the cost of moving a unit of mass from a point $x \in X$ to a point $y \in Y$ and the integral 
\[\int_{X \times Y} c(x,y) \, d \pi(x,y) \]
is then the total cost of transportation according to the measure $\pi$. Which cost function one uses depends on the context of the problem, but is fixed as part of the problem statement. In Monge's work, $X$ and $Y$ were taken to be sets in Euclidean space and the cost function was simply the distance $c(x,y) =|x-y|$. For many applications however, it is preferable to instead use the squared-distance $c(x,y) =|x-y|^2$. 

The Kantorovich framework is a flexible and powerful approach to study optimal transport for several reasons. First, an \textit{optimal coupling} exists under very general circumstances.

\begin{theorem} \label{Existence}
Suppose that $\psi : X \to \mathbb{R} \cup - \infty$ and $\phi : Y \to \mathbb{R} \cup - \infty$ are two upper-semicontinuous functions with $\psi \in L^1(\mu)$ and $\phi \in L^1(\nu)$, such that $c(x,y) \geq \psi(x) + \phi(y)$. Then there is a coupling $\pi \in \Pi(\mu, \nu)$ which minimizes the total cost.
\end{theorem}
In particular, if the cost function is non-negative, an optimal coupling exists (in this case, we take $\psi$ and $\phi$ to be identically zero).

Second, the Kantorovich problem has a dual formulation, which plays an essential role in proving both theoretic and numeric results in optimal transport. Take $\tilde{c}(x,y)=-c(x,y)$, and we consider the problem of maximizing 
\begin{eqnarray*}
   \sup_{\pi \in \Pi(\mu,\nu)} \int_{X \times Y}  \tilde{c}(x,y) \, d \pi(x,y) =
\end{eqnarray*}
\begin{eqnarray*}
    \inf_{\widetilde{\phi}, \widetilde{\psi}} \left \{   \left( \int_X \widetilde{\psi}(x) d \mu(x) +  \int_Y \widetilde{\phi}(y) d \nu(y) \right)  ~ \Big{ \vert} ~ \begin{aligned} \widetilde{\psi} \in L^1(\mu), ~ \widetilde{\phi} \in L^1(\nu),  \\ \widetilde{\psi}(x) + \widetilde{\phi}(y) \geq \tilde{c}(x,y)
    \end{aligned}\right \} .
\end{eqnarray*}  
An economic interpretation the Kantorovich dual problem is clearly explained in Villani's textbook (\cite{villani2009optimal}, Chapter 5). 
Under relatively mild assumptions (for instance, when $c$ is continuous, real-valued and satisfies some mild integrability assumptions), the infimum on the right hand side is obtained. 

Treating $\tilde{c}\equiv -c$ as the gain, the maximizing and minimizing versions of the optimal transport problem are equivalent,  with the optimal couplings and solutions unchanged (apart from $\widetilde{\psi} \equiv -\psi, \widetilde{\phi} \equiv - \phi$). This maximizing notation is more congruent with information geometric convention of using convex functions and treating the Legendre transform as a maximization (instead of minimization) problem. For notational simplicity, we will drop the $\, \widetilde{\mbox{}} \,$ on the functions $\psi, \phi$ in the rest of the paper.  See the Remark in the next subsection.

\subsection{Solving the Kantorovich dual problem}

Solutions to the Kantorovich dual problem are linked to a duality which appears throughout optimal transport, dubbed {\it $c$-duality}.  More precisely, the functions $\psi$ and $\phi$ in the formulation of Kantorovich problem are linked: $\phi = \psi^c$, where $\psi^c$ is the \textit{$c$-conjugate} of $\psi$, defined as 
\begin{equation}
\label{c-conjugatedef}
 \psi^c(y) = \sup_{x \in X} \left(  \tilde{c}(x,y) - \psi(x) \right) 
\end{equation}
or equivalently,
\begin{equation}
\label{c-conjugatedef2}
\psi^c(y) = - \inf_{x \in X} \left( c(x,y) + \psi(x)  \right). 
\end{equation}
Furthermore, the optimal coupling $\pi$ is concentrated on the $c$-subdifferential of $\psi$, which is defined as follows:
\[ \partial_c \psi = \left \{ (x,y) \in X \times Y \,|\, \psi^c(y)+\psi(x) = \tilde{c}(x,y) \right \} . \]

Given an arbitrary $\psi$, the function $\psi^c$ defined by \eqref{c-conjugatedef} or 
\eqref{c-conjugatedef2}
is said to be $c$-convex, and a $c$-convex function is exactly one that can be realized as the $c$-conjugate of some other function.

Taking the $c$-conjugate again, we find that
\begin{equation*}
 \psi^{cc}(x) = \sup_{y \in Y} \left(  \tilde{c}(x,y) - \psi^c(y) \right) 
\end{equation*}
or 
\begin{equation*}
\psi^{cc}(x) = - \inf_{y \in Y} \left( c(x,y) + \psi^c(y)  \right).
\end{equation*}
This pattern continues for further applications of $c$-conjugation. It is readily shown (see the tutorial by Ambriosio and Gigli \cite{ambrosio2013user}) that
\begin{enumerate}
    \item[(i).] $\psi^{ccc} = \psi^{c}$;
    \item[(ii).] $\psi^c$ is always $c$-convex;
    \item[(iii).] $\psi^{cc}=\psi$ when $\psi$ is $c$-convex. 
\end{enumerate}

The $c$-conjugate transform is a generalization of Legendre duality: for the cost $c(x,y)=-x \cdot y$, that is, $\tilde{c}(x,y) = x \cdot y$, $c$-duality becomes the regular Legendre duality: the $c$-subdifferential is the usual subdifferential for a convex function, and the $c$-conjugate is precisely the standard Legendre transformation. This link between optimal transport and the analysis of $c$-convex functions lies at the heart of optimal transport. 

\subsubsection*{Remark} \label{Remark} The language of $c$-duality, as used in optimal transport, is often cast in $c$-concave functions by standard references. In information geometry which invokes the regular Legendre duality, convex functions are adopted for convenience. Of course, a convex function can be made concave simply by multiplying it by $(-1)$, i.e., changing its sign; however this is not the case of $c$-convexity and $c$-concavity of functions. To facilitate comparison in this article, we restate the results of $c$-duality in terms of $c$-convexity instead of $c$-concavity. 

\subsection{Existence theory to the Monge problem}
\label{Existence of Monge problem}
In Monge's work, he made the additional assumption that $\pi$ is concentrated on the graph of a measurable map $T:X \to Y$. In other words, each point in $X$ is sent to a unique point in $Y$.
A solution to the Monge problem need not exist in general. For instance, if $\mu$ is atomic but $\nu$ is continuous, then there are no transport maps at all, let alone an optimal one.

One great success of the modern theory of optimal transport is to develop a satisfactory existence theory for the Monge problem. In particular, Brenier \cite{brenier1987decomposition}, Gangbo-McCann \cite{gangbo1996geometry}, and Knott-Smith \cite{knott1984optimal} established that for fairly general cost functions and measures, the solution to the Kantorovich problem is actually a solution to the Monge problem.

\begin{theorem}\label{Deterministictransport} 
Let $X$ and $Y$ be two open domains of $\R^n$ and consider a cost function $c:X \times Y \to \R$. Suppose that $ d \mu$ is a smooth probability density with respect to the Lebesgue measure supported on $X$ and that $d \nu $ is a smooth probability density with respect to the Lebesgue measure supported on $Y$. Suppose that the following conditions hold:

\begin{enumerate}
\item[i.] the cost function $c$ is of class $C^4$ with $\| c \|_{C^4(X \times Y)} < \infty$;
\item[ii.] for any $x \in X$, the map $ Y \ni y \to  c_x(x,y) \in \R^n$ is injective;
\item[iii.] for any $y \in Y$, the map $ X \ni x \to  c_y(x,y) \in \R^n$ is injective\footnote{Items ii. and iii. are often known collectively as the \emph{Twist} condition.};
\item[iv.] $\det(c_{x,y})(x,y) \neq 0$ for all $(x,y) \in X \times Y$.
\end{enumerate}
Here the subscript of $c$ denotes taking partial derivative, with respect to $x$ if before the comma and with respect to $y$ if after the comma. 

Then, there exists a $c$-convex function $\psi: X \to \R$ such that the map $\mathbb{T}_\psi : X \to Y$ defined by $\mathbb{T}_\psi(x) := c\textrm{-}\exp_x(\mathsf{D} \psi(x))$ is the unique optimal transport map sending $\mu$ onto $\nu$. Furthermore, $\mathbb{T}_\psi$ is injective $d \mu$-a.e. and is an Alexandrov solution to the Jacobian equation
\begin{equation} \label{Monge Ampere}
 | \det(\mathsf{D} \mathbb{T}_\psi(x))| = \frac{ d \mu (x)}{ d \nu(\mathbb{T}_\psi(x))} \hspace{.2in} d \mu-a.e.
\end{equation}
\end{theorem}

In order to express Equation \eqref{Monge Ampere} more concretely, we recall the notion of the $c$-exponential map (denoted $c\textrm{-}\exp_x$).

\begin{definition} [$c$-exponential map]
For any $x \in X, y \in Y, p \in \mathbb{R}^n$, the $c$-exponential map satisfies the following identity.
\[ c\textrm{-}\exp_x(p) = y \iff p = - c_x(x, y) = \tilde{c}_x(x,y). \]
\end{definition}

The name $c$-exponential comes from the fact that on a Riemannian manifold with the cost $c=d^2(x,y),$ this coincides with the standard exponential map.  For a more complete reference on the existence of solutions to the Monge problem, see Chapters 10-12 of \cite{villani2009optimal}.

\subsection{The regularity of the Monge map \& the MTW tensor}
\label{sec_regularity}

As shown in Theorem \ref{Deterministictransport}, for fairly general costs and measures, the solution to the Kantorovich problem will also be a solution to the Monge problem, and it is induced by the $c$-subdifferential of a potential function $\psi$. Furthermore, this potential function will be a weak solution to a fully non-linear degenerate-elliptic PDE known as the {\it Jacobian equation}. From here, the natural follow-up question is to understand the behavior of the solution --- is it continuous? injective? smooth?

To answer these questions, one must find a priori estimates for the Jacobian equation. A full discussion of this theory is outside the scope of this survey, but let us provide a brief summary. 
 A priori, the potential function $\psi$ will be $c$-convex
 and locally Lipschitz, which implies that its $c$-subdifferential (and thus the transport map) is defined Lebesgue almost-everywhere. However, the transport map need not be continuous, even for smooth measures. For example, for the squared-distance cost $c(x,y)= |x-y|^2$,  Caffarelli \cite{caffarelli1992regularity} showed that for the transport between arbitrary smooth measures to be continuous, it is necessary to assume that the target domain $Y$ is convex. In other words, non-convexity of $Y$ is a \textit{global obstruction} to regularity.\footnote{There are easier ways to see that the optimal transport map might fail to be continuous. For instance, if $X$ is connected whereas $Y$ is disconnected (and the measures have full support), there are no continuous transport maps at all, let alone an optimal one.}

From an analytic perspective, the reason for this somewhat pathological behavior is that Equation \eqref{Monge Ampere} is \textit{degenerate-elliptic}. To explain its meaning, let us specialize to the squared-distance cost in Euclidean space, in which case the $c$-convexity of $\psi$ is convexity in the usual sense and the Jacobian equation simplifies to the following:
\begin{equation} \label{MongeAmpere2} \det \left( \mathsf{D}^2 \psi(x) \right) = \frac{d \mu(x)}{d \nu(\mathsf{D} \psi (x))}. \end{equation}

For a convex function $\psi$, the linearized Monge-Amp\'ere operator
\[ Lf = \sum_{i,j}  \left( \left(\operatorname{det} \mathsf{D}^{2} \psi\right)\left(\mathsf{D}^{2} \psi\right)^{-1} \right )^{ij} \frac{\partial^2} {\partial x^i \partial x^j} f  \]
 is non-negative definite. Furthermore, the potential $\psi$ will satisfy Equation \eqref{MongeAmpere2} where we replace the left-hand side by its linearized Monge-Amp\'ere equation $L$. As such, to establish regularity for $\psi$, we need to derive estimates for solutions to
 \[ L\psi = \frac{d \mu(x)}{d \nu(\mathsf{D} \psi (x))}, \]
which is a quasi-linear equation.
 
 However, without an {\it a priori} bound on $\mathsf{D}^2 \psi$, $L$ can have arbitrarily small eigenvalues, which prevents the use of standard elliptic theory to prove regularity of $\psi$. On the other hand, if we were able to establish an {\it a priori} $C^2$ estimate on $\psi$, then the linearized operator would be uniformly elliptic, and we would be able to use Schauder estimates to show that the transport is smooth.

This phenomena gives rise to a striking dichotomy. When the cost functions and measures satisfy the hypotheses of Theorem \ref{Deterministictransport} and the measures are smooth (i.e., $C^\infty$), the optimal transport map is either $C^\infty$-smooth, or else discontinuous.\footnote{This dichotomy only holds in the interior. It is possible for the Monge map to be smooth in the interior, but become singular near the boundary.} In the latter case, work of Figalli and De Phillippis \cite{de2015partial} shows that the discontinuity occurs on a set of measure zero and that the transport is smooth elsewhere. Kim and Kitagawa  \cite{kim2016prohibiting} showed that for costs which satisfy the MTW(0) condition (Definition \ref{MTWkappadefinition}), there cannot be isolated singularities. However, there are many open questions about the structure of the singular set.

\begin{question}[De Phillippis-Figalli \cite{de2014monge}]
For an arbitrary smooth domain and cost function, does the singular set for the transport have Hausdorff dimension at most $n-1$? If so, is this set rectifiable?
\end{question}

\subsubsection{Local obstruction to regularity}

In the 1990s, Caffarelli, Delanoe, Urbas and several others developed a theory for the regularity of the Monge problem when the cost function is the squared distance in Euclidean space (see, e.g., \cite{caffarelli1992regularity,delanoe1991classical,urbas1997second}).
Though imposing the additional assumption that the target space $Y$ was convex, this line of work was able to apply techniques from Monge-Amp\'ere equations, due to its simplified form \eqref{MongeAmpere2}, to derive the relevant {\it a priori} estimates. However, for more general cost functions, the problem of regularity remained open. 

In 2005, a breakthrough paper by Ma, Trudinger, and Wang  \cite{ma2005regularity} established smoothness for the transport for general costs under two additional assumptions. The first was a global condition, that the supports of the initial and target measures are \textit{relatively c-convex}. The second was a local condition, that a certain fourth-order quantity, known as the \textit{MTW tensor}, is positive.  Later, Trudinger and Wang \cite{trudinger2009second} extended their previous work to the case when this tensor is non-negative (i.e., when $\mathfrak{S}(\xi,\eta)\geq 0$). 

\begin{definition}[MTW($\kappa$) condition] \label{MTWkappadefinition}
A cost function $c$ is said to satisfy the MTW($\kappa$) condition if for all vector-covector pairs $(\xi, \eta)$ with $ \eta(\xi)=0$, the following inequality holds.
\begin{equation}\label{MTW0}
\mathfrak{S}(\xi,\eta) := \sum_{i,j,k,l,p,q,r,s} (c_{ij,p}c^{p,q}c_{q,rs}-c_{ij,rs})c^{r,k}c^{s,l} \xi^i \xi^j \eta^k \eta^l \geq \kappa \|\eta\|^2 \| \xi \|^2 
\end{equation}
Here, the notation $c_{I,J}$ denotes 
$\partial_{x^I} \partial_{y^J} c$ for multi-indices $I$ and $J$. Furthermore, $c^{i,j}$ denotes the matrix inverse of the mixed derivative $c_{i,j}$. Finally, all the derivatives and the MTW tensor are understood to be taken at the point $(x,y)$.
\end{definition}

It is worth remarking that in \eqref{MTW0}, all of the terms involve at least a third derivative, so $\mathfrak{S}$ identically vanishes for the squared-distance cost in Euclidean space, which explains why it did not appear in the earlier work on the regularity problem.

Although the Ma-Trudinger-Wang theory was a significant breakthrough, it also raised many questions. For instance, the geometric significance of $\mathfrak{S}$ was not well understood and it was unclear whether it plays an essential role in optimal transport or if it was merely a technical assumption in Ma-Trudinger-Wang's work. 

These questions were studied by Loeper \cite{loeper2009regularity}, who showed that for costs which are $C^4$, the MTW(0) condition is equivalent to requiring that the $c$-subdifferential of an arbitrary $c$-convex function be connected. This provides a geometric interpretation of the MTW tensor in terms of convex analysis, and shows that the failure of MTW condition is a \textit{local obstruction} to regularity. More precisely, this work shows for cost functions which do not satisfy the MTW condition, the space of smooth optimal transports fails to be dense within the space of all optimal transports. In other words, given a cost which does not satisfy the MTW(0) condition, it is possible to find smooth measures $\mu$ and $\nu$ for which the Monge optimal transport is discontinuous.

There are many open questions remaining in the regularity theory of optimal transport, and we will mention several in Section \ref{Geometrizing Optimal Transport section}. However, let us mention one prominent conjecture due to Villani.

\begin{conjecture}[Villani's Conjecture]
If $(M,g)$ is a Riemannian manifold such that the cost $c(x,y)=d(x,y)^2$ satisfies the MTW condition, then all of the injectivity domains of $M$ are convex.
\end{conjecture}

\subsection{The Wasserstein metric}

Though optimal transport is well defined for a very general class of cost functions, many important applications consider the case when the cost is some power of the distance between points on a metric space $(X,d)$. This gives rise to the notion of the \textit{Wasserstein metric}, which provides a notion of distance between probability measures. The name ``Wasserstein metric" is also referred as the \textit{Kantorovich-Rubinstein distance}, the \textit{Monge-Kantorovich distance} (both of which might be more historically accurate) or the {\it Earth Mover Distance} (which is much more descriptive) \cite{santambrogio2015optimal}. 
 
To define the Wasserstein metric, we consider the space of probability measures with finite $p$-th moment:
\begin{equation} \label{Wasserstein domain}
    \mathbb{P}_p(X) = \left \{ \textrm{probability measures} \, \mu \,  \Big{ \vert} \, \int_X d(x,x_0)^p d \mu(x) < \infty \textrm{ for any } x_0 \in X \right \} 
\end{equation}
For two probability measures $\mu, \nu \in \mathbb{P}_p(X)$, we can define their $p$-Wasserstein distance as 
\[ W_p(\mu, \nu) = \left( \inf_{\pi \in \Pi(\mu,\nu)} \int_{M \times M} d(x,y)^p d \pi \right)^{1/p} . \]
This is a genuine metric in that it is positive-definite, finite, symmetric, and satisfies the triangle inequality. Note that here ``metric'' is used in the sense of a distance function, and is different from its usage of Riemannian metric (an inner product) or a Finsler metric (a smooth, homogeneous and subadditive ``norm'' function) at each point of the manifold.

\subsection{Displacement interpolation}

Displacement interpolation is a continuous-time version of optimal transport which interpolates between the initial and final measure \cite{mccann1997convexity}. For this, one must assume that the source space $X$ and the target space $Y$ are the same. However, the setting for displacement interpolation is quite general otherwise, and can be defined when $X$ is a Polish metric space and the cost functions is of the form
\begin{equation}
c(x, y)=\inf \left\{\mathcal{L}(\gamma) \, | \, \gamma_{0}=x, \,\, \gamma_{1}=y  \right\}
\end{equation}
where $\gamma$ is a path connecting $x$ to $y$ and $\mathcal{L}$ is a Lagrangian action
\begin{equation*}
\mathcal{L}(\gamma)=\int_{0}^{1} L\left(\gamma_{t}, \dot{\gamma}_{t}, t\right) d t.
\end{equation*}

In this section, we will focus on the case where the cost is the distance to some power, so will specialize to the associated action
\begin{equation*}
\mathcal{L}_p(\gamma)=\int_{0}^{1} \frac{1}{p} \, \left \vert \dot \gamma_{t} \right \vert^p  d t.
\end{equation*}

Using the action $\mathcal{L}_p$, it is possible to induce on $\mathbb{P}_p(X)$ (where $(X,d)$ is {\it a priori} merely a metric space) the structure of a \emph{length space}, where the distance $d(x,y)$ can be realized as the infimum as the lengths of paths connecting $x$ and $y$ as end-points:
\[ d(x,y) = \inf_{\gamma} \left \{ \ell(\gamma) ~\big \vert \, \gamma_0 = x, \, \gamma_1 = y \right\} \]
where $\ell(\gamma)$ is defined to be the intrinsic length of the path
\[ \ell(\gamma) =  \sup_{P} \left \{ \sum_{P}  d\left(\gamma_{p_{i-1}}, \gamma_{p_{i}} \right ) ~|~ P \textrm{ a partition of $[0,1]$ with } p_i \in P \right \}.
\]

\begin{theorem}[\cite{villani2009optimal} Corollary 7.22] Let $(X, d)$ be a complete separable, locally compact length space. Then, given any two measures $\mu_{0}, \mu_{1} \in \mathbb{P}_{p}(X)$, and a continuous curve $\left(\mu_{t}\right)_{0 \leq t \leq 1}$ in $\mathbb{P}_p(X)$, the following properties are equivalent:
\begin{enumerate}
    \item   $\mu_{t}$ is the pushforward by the flow of  minimizing, constant-speed geodesics from $x$ to $\mathbb{T}(x)$, where $\mathbb{T}:X \to X$ is the solution to the optimal transport problem.\footnote{If the solution to this transport problem is not Monge, this should be interpreted in the sense of couplings where mass may split at the initial time along various geodesics.}
   \item $\left(\mu_{t}\right)_{0 \leq t \leq 1}$ is a geodesic curve in the space $\mathbb{P}_{p}(X)$.
\end{enumerate}
Moreover, if $\mu_{0}$ and $\mu_{1}$ are given, there exists at least one such curve.
\end{theorem}

In other words, displacement interpolation provides geodesics in the space of probability measure, and this interpolation is induced by geodesics in the underlying space $X$. Furthermore, one can understand the geometry of $X$ by studying the behavior of the measures along displacement interpolation. For instance,  the convexity of certain entropy functionals along geodesics in $\mathbb{P}_{p}(X)$ is equivalent to Ricci lower bounds\footnote{More precisely, most of these bounds involve both the curvature and the dimension, so are known as curvature-dimension bounds.} on the original space $X$ \cite{von2005transport,ohta2011displacement}.

\subsection{Otto calculus}
\label{otto_cal}
While the length space structure of $\mathbb{P}_p(X)$ holds for arbitrary values of $p$, the case $p=2$ is distinguished. When the space $X$ is a Riemannian manifold, Otto \cite{otto2001geometry} showed that the space $\mathbb{P}_2(X)$ (metrized by the $2$-Wasserstein distance) also admits the structure of a formal Riemannian metric. 
For ease of exposition, let us look at the special case of $X \subseteq \mathbb{R}^n$ a part of Euclidean space, and consider the subset of $\mathbb{P}_2(\mathbb{R}^n)$ which can be represented as probability density functions over $\mathbb{R}^n$ with respect to the Lebesgue measure $dx\equiv d \mbox{vol}_x$:
 \[ \mathcal{M} = \mathbb{P}_2(X)=\left\{ \rho \, | \, \rho(x) > 0 , \, \forall  x \in X \subset \mathbb{R}^{n}, \, \int_{X} \rho \, dx =1 \right\}.  \]
More precisely, the function $\rho$ is the Radon-Nikodym derivative of a measure which is Lebesgue absolutely continuous. 
One can then consider the tangent space at a point $\rho \in \mathbb{P}_2(X)$ as
 \[ T_{\rho} \mathcal{M}=\left\{s \, | \, s \text { is a function on } X, \, \int s dx =0 \right\} .\]
 To see why this can be understood as the tangent space, consider the ``nearby" probability measure $\rho_\epsilon = \rho + \epsilon s$. The integral of $\rho_{\epsilon}$ will be 1: $\int \rho_\epsilon \, dx = 1$. Furthermore, by restricting $s$ and $\rho$ to reside in a suitable function space (e.g., $C_c^\infty(\mathbb{R}^n)$) and to satisfy $\textrm{supp} (s) \subset \textrm{supp}(\rho)$, we can ensure that for $\epsilon$ sufficiently small, $\rho + \epsilon s$ is non-negative. So 
 $$
 s(x) = \lim_{\epsilon \rightarrow 0} \frac{\rho_\epsilon(x) - \rho(x)}{\epsilon}.
$$

One can then use the fluid-dynamic elliptic equation 
\begin{equation} \label{elliptic}
-\mathsf{div} \cdot(\rho \, \mathsf{grad} \, \phi)=s 
\end{equation}
to obtain a second set of ``coordinates" for the tangent fibers, which is denoted $\phi$:
\[ T_{\rho} \mathcal{M} \cong \left\{\text { functions } \phi \text { on } \mathbb{R}^{n}\right\} / \{ \phi_1 -\phi_2 \equiv const \} .\]

In other words, functions which differ by constants are deemed identical and are then identified with the original coordinates $s$ via the solution to an elliptic equation. Here, since this computation is merely formal, one omits any discussion of which Sobolev spaces each of these families of functions reside in, or what the boundary conditions should be. 

In this setting, Otto was able to prescribe a Riemannian metric $g$ on $\mathbb{P}_2(X)$, defined as the quantity
 \[ g_{\rho}\left(s_{1}, s_{2}\right)=\int_{X}  (\mathsf{grad} \, \phi_{1} \cdot \mathsf{grad} \, \phi_{2}) \, \rho \, dx \,\, . \]
Here, the dot-product is the standard inner-product operation on $X=\mathbb{R}^n$. The Riemannian metric $g$ induces the 2-Wasserstein distance on $\mathbb{P}_2(\mathbb{R}^n)$ as its distance function, so that $\mathbb{P}_2(\mathbb{R}^n)$ admits a formal Riemannian structure. If $X$ is a Riemannian manifold instead of a subset of Euclidean space $\mathbb{R}^n$, then the dot-product $\cdot$ can be replaced by a Riemannian metric $\langle \cdot, \cdot \rangle$ {\it on $X$}, and Equation \eqref{elliptic} is defined in terms of this metric (on $X$) as well.  Before moving on, it is worth mentioning contemporaneous work of Benamou and Brenier \cite{benamou2000computational}, which studied this structure from a fluid mechanical perspective.


One ought to clearly distinguish the two spaces involved in Otto's framework: 
\begin{itemize}
\item $X$, which we call the ``downstairs'' manifold, that is a Riemannian manifold itself endowed with a Riemannian metric $\langle \cdot, \cdot \rangle$, a special case being a Euclidean space; each point $x\in X$ is the instantiation of a random variable;
\item $\mathbb{P}_2(X)$, which we call the ``upstairs'' space and which admits a formal Riemannian metric $g_\rho(\cdot, \cdot)$; each point $\rho \in \mathcal{M} = \mathbb{P}_2(X)$ is a probability measure on $X$.  
\end{itemize}
See Section \ref{Upstairs and downstairs section} for examples. So
\begin{itemize}
\item $\mathsf{grad}$ and $\mathsf{div}$ are operators on the downstairs manifold $X$;
\item $s$ and $\phi$ are functions on $X$, linked through Equation \eqref{elliptic} with respect to any given $\rho$ (i.e., treating $\rho$ as fixed);
\item $\rho$ is an element of $\mathcal{M} = \mathbb{P}_2(X)$ while $s$ is an element of $T\mathcal{M}$. 
\end{itemize}

Otto calculus plays a central role in modern optimal transport --- one of its hallmarks is to show that Ricci curvature lower bounds of $X$ gives rise to geodesic convexity of certain entropy functionals on $\mathbb{P}_2(X)$. Using ideas from Otto calculus, Lott, Sturm, Villani and others have developed notions of Ricci curvature bounds on non-smooth metric-measure spaces. For an overview of this theory, we refer the reader to the following survey \cite{villani2016synthetic}. We also refer to the following paper of Lott \cite{lott2008some} for curvature calculations of the Wasserstein metric. It is worth emphasizing that the Riemannian structure is merely formal, as $\mathbb{P}_2(X)$ is generally not a Banach manifold. As such, many of the rigorous proofs in this line of work do not directly use the Otto calculus, but instead follow a motto of Villani to ``[t]hink Eulerian, prove Lagrangian" (\cite{villani2009optimal}, page 444). In other words, one uses the Otto calculus to derive results, but establishes them rigorously using displacement interpolation.

In the case where the ``upstairs'' space of probability measures $\mathbb{P}_2(X)$ admits a Kullback-Leibler divergence, then the gradient flow (on the upstairs space, hence subscripted by $W_2$) of this KL divergence can be calculated:
$$
\frac{\partial \rho_t}{\partial t} = - \mathsf{grad}_{W_2} K(\rho_t||\rho_0) ,
$$
which gives rise to the Fokker-Planck equation 
$$
\frac{\partial \rho_t}{\partial t} =
\nabla \cdot (\rho_t \nabla \phi) + \Delta \rho_t .
$$
where $\rho_0 \simeq e^{-\phi}$ is the Gibbs measure (see Chapter 24 of \cite{villani2009optimal}, and see \cite{jordan1998variational} for a proof of convergence for this equation). 
There are many results relating well-known PDEs to the gradient flows of certain functionals in $\mathbb{P}_2(X)$ (or more generally $\mathbb{P}_p(X)$ \cite{otto2001geometry}), and this is a thriving area of functional analysis.

\section{Linking Optimal Transport to Information Geometry}
In this section, we survey some areas of active research at the interface of information geometry and optimal transport. Due to the anticipated familiarity with information geometry by the audience of this Journal, we do not provide any review of the framework of information geometry (see references like \cite{amari2012differential,AN00,ay2017information}). Rather, we start by addressing some key differences between Wasserstein geometry and information geometry.  

\subsection{Comparing Wasserstein geometry with information geometry} 
Both Wasserstein geometry and information geometry provide ways to quantify the proximity of probability measures, the former through solutions to optimal transport and the latter using divergence (or contrast) functions, such as Kullback-Leibler (KL) divergence (or more generally, Csiszar's $f$-divergence \cite{csiszar1964informationstheoretische}). A divergence function is not a distance, in that it is not necessarily symmetric or satisfy the triangular inequality. However, it is still a non-negative function of two probability density functions  $p$ and $q$ (or probability measures $\mu$ and $\nu$), and equals zero if and only if $p=q$ almost surely.

There are several important distinctions between Wasserstein distance and divergence functions. Firstly, they are defined on different classes of probability measures. The KL divergence is only defined for probability measures $\mu, \nu$ are which are absolutely continuous with respect to some reference measure $\rho_0$ (or at least one with respect to another). More precisely, on a sample space $X$, we can express the KL divergence as the following
\begin{equation}\label{KL divergence integral}
 K[p||q] = \int_X \frac{d\mu}{d\rho_0} \log \left(\frac{d\mu}{d\rho_0} / \frac{d\nu}{d\rho_0} \right) \, d\rho_0 = \int_X \log \frac{d\mu}{d\nu} \, d\mu .
\end{equation}
As such, it is necessary to restrict our attention to probability measures whose Radon-Nikodym derivatives (with respect to $\rho_0$) exist and live in an appropriate Orlicz space (so that the integral \eqref{KL divergence integral} converges). 

On the other hand, the $p$-Wasserstein distance is defined for any two probability measures with finite $p$-th moment ($1\le p < \infty)$, without any assumptions about absolute continuity. In fact, the source and target distributions may well have different supports on the sample space $X$. In other words, the Wasserstein distance can be finite even when the support of $\mu$ and $\nu$ is disjoint, which is not the case for the KL divergence.  


\subsubsection{``Vertical'' versus ``horizontal'' displacement}

Transport-based Wasserstein distances and information geometric divergence functions both measure the ``distance'' between two probability measures/densities, but in very different ways.

\begin{figure}
    \centering
    \includegraphics[width=.9\linewidth]{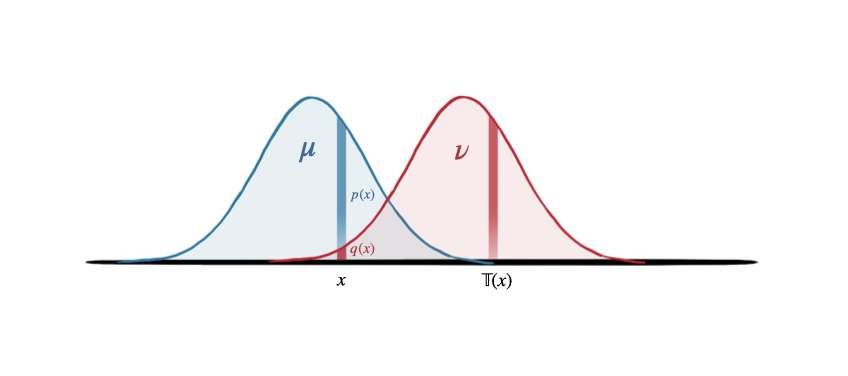}
    \caption{Vertical versus horizontal displacement}
    \label{fig:verticalandhorizontaldistances}
\end{figure}

To explain this, let us first reframe the notion of a divergence in a slightly different way. We suppose that $\mu$ and $\nu$ are absolutely continuous with respect to a reference measure $\rho_0$ with Radon-Nikodym derivatives $p(x) = \frac{d \mu}{d \rho_0}$ and $q(x)= \frac{d \nu}{d \rho_0}$.
In this setting, we can express the divergence $D$ between $\mu$ and $\nu$ in terms of the following integral:
$$
D[p||q] = \int k(p(x), q(x)) \, d\rho_0(x) 
$$
where $k(a,b)$ is some non-negative smooth function on $\R_{\geq 0}^2$ which is equal to zero if and only if $a=b$. For instance, to construct the KL-divergence, we use the function
$$ 
k(p(x), q(x)) = q(x) - p(x) - p(x) \log \frac{q(x)}{p(x)},  
$$
which satisfies $k(p(x), q(x)) \geq 0 , \, \forall x \in X .$ Divergence functions of this form are sometimes called {\it extended KL divergences} \cite{zhu1995information}, and their importance for generating the dualistic information geometry was recognized and made use of extensively by the second named author \cite{zhang2004divergence}. Furthermore, the notion of gauge freedom \cite{Naudts2018Rho} can be introduced, allowing an arbitrary monotone function to be used in defining the divergence function. 

From this integral representation, we see that divergences $D$ compare the ``heights'' $p(x)$ and $q(x)$ of the probability density functions at each point $x$, and use the non-negative function $k(p(x), q(x))$ to measure the difference in these heights. Then, we integrate these differences over the random variable $x$ with respect to the reference measure $\rho_0$. In this way, divergences $D$ measure displacement accumulated ``vertically.'' From this, it immediately follows that these divergences remain invariant under measure-preserving mappings of the underlying sample space (for example, when the random variable $x$ in the sample space is relabeled, i.e., $X$ undergoes a permutation). 

On the other hand, the Wasserstein distance is defined by the cost of optimally transporting one of the measures to the other. As such, this is ``horizontal" measure of displacement, in which mass from the source measure $\mu$ is moved horizontally through space to the target measure $\nu$. Since conjugation by a measure-preserving mapping generally does not preserve the cost function or the optimal map, measure-preserving maps will not preserve Wasserstein distances. Furthermore, on a Riemannian manifold one can find measure-preserving mappings whose conjugation make the Wasserstein distance arbitrarily small. 

\subsubsection{Topological considerations}
\label{Differing topologies}
Given that they are defined for different classes of probability measures, it might be expected that the $p$-Wasserstein distance and KL divergence induce different topologies, and this is indeed the case. Due to the asymmetry of the KL divergence, there are actually two separate topologies depending on whether one considers open balls with respect to the first or second argument (see, e.g., \cite{belavkin2015asymmetric}). 
However, convergence in either topology implies convergence in total variation, so both are fairly strong notions of convergence which are invariant under measure-preserving mappings.

On the other hand, the Wasserstein distance metrizes the weak topology of $\mathbb{P}_p$, in that convergence in Wasserstein distance implies weak convergence of measures and convergence of their first $p$-th moments (Theorem 6.9 \cite{villani2009optimal}).

As a result of this, we generally (but not always) expect that the topology induced by the Wasserstein distance will be coarser than that induced by the Kullback-Liebler divergence. As such, we can often find bounds on the Wasserstein distance in terms of the KL divergence, but essentially never the other way round (except for finite spaces).  However, when the measures are further controlled in a stronger sense, it is sometimes possible to prove interpolation results between the Wasserstein distance and KL divergence. Most famously, this idea shows up in the ``HWI-inequalities"  of Otto and Villani \cite{otto2000generalization}. 

\subsection{Entropy-regularized transport}
\label{Entropy-regularized transport}
One link between optimal transport and information geometry comes in the form of entropy-regularized optimal transport. 
To explain this, it will be convenient to specialize to the case of discrete probability measures (instead of probability density functions) which are supported on finite sets.\footnote{See Section \ref{Works in Otto calculus/IG} for related work in the continuous setting.} Assuming the source is supported on $X=\{1, \cdots, n\}$ and the target on $Y=\{1, \cdots, m\}$, with a cost $C_{ij}$ incurred from transferring from $i \in X$ to $j \in Y$. The coupling can be written as
$$
\Pi(p,q) = \left \{ P_{ij}: P_{ij} \ge 0, \, \sum_{ij} P_{ij} = 1, \, \sum_j P_{ij} = p_i, \, \sum_i P_{ij} = q_j \right \} . 
$$
Here, 
$$
\sum_i p_i = \sum_j q_j = 1 ,
$$
and the objective function for transport is:
$$
W(p,q) = \inf_{P \in \Pi(p,q)} \langle P, C \rangle = \inf_{P \in \Pi(p,q)} \sum_{i,j} P_{ij} C_{ij} .
$$
As a linear-programming problem, there are standard method for solving optimal transport over discrete probabilities, but this may be computationally inefficient.

Cuturi \cite{cuturi2013sinkhorn} proposed a relaxation to the above problem by introducing a regularization term using the Shannon entropy $H(P)$:
\begin{eqnarray} 
W_\epsilon(p,q) &=& \inf_{P \in \Pi(p,q)} \langle P, C \rangle - \epsilon H(P) \nonumber \\
&=& \inf_{P \in \Pi(p,q)} \sum_{i,j} P_{ij} C_{ij} + \epsilon \sum_{i,j} P_{ij} \log P_{ij} .
\label{cuturi}
\end{eqnarray}

Applying the Lagrange multiplier method to this entropy-regularized (also called entropy-relaxed) problem, Cuturi showed that the solution is of the form
\begin{equation} \label{opt_sol}
P^{opt}_{ij} = \exp \left [  \frac{1}{\epsilon} \left(\alpha_i + \beta_j - C_{ij} \right) \right] = u_i K_{ij} v_j 
\end{equation}
where $\alpha$ and $\beta$ are vectors of Lagrange multipliers, which are specified by the following equations:
$$
\sum_j P^{opt}_{ij} = p_i, \,\,\,\, \sum_i P^{opt}_{ij} = q_j , \quad \forall i, j
$$
and the matrix $K$ and vectors $u$ and $v$ denote
$$
K_{ij} = \exp \left(-\frac{1}{\epsilon} C_{ij}  \right) , \,\,\,\,
u_i = \exp\left(\frac{\alpha_i}{\epsilon} \right), \,\,\,\,
v_j = \exp \left( \frac{\beta_j}{ \epsilon} \right) .
$$
Note that using entropy $H(P)$ as the regularization term is equivalent to using KL divergence $K[P||p \otimes q]$ as a regularizer. Though $W_\epsilon(p,q)$ may no longer be a distance, it approaches the Wasserstein distance $W(p,q)$ in the limiting case of $\epsilon$ tending zero. 

Amari, Karakida, Oizumi \cite{amari2018information} observed that the solutions $P^{opt}$ from Equation \eqref{opt_sol} form an exponential family, where $p,q$ are the expectation parameters, and $\alpha, \beta$ are the natural parameters. Viewed in this way, the Sinkhorn algorithm used to iteratively solve \eqref{opt_sol} has the information geometric interpretation of successive applications of the $e$-projection to a pair of $m$-autoparallel submanifolds. One can also define a Bregman-like divergence function on these autoparallel submanifolds (the so-called Cuturi function).
A similar geometric interpretation, in terms iterative Bregman projections, was given by Benamou et al. \cite{benamou2015iterative}. 

Muzellec et al.\ \cite{muzellec2017tsallis} extended Cuturi's entropy-regularized framework to the optimal transport problem by using Tsallis entropy. Of course, more general regularizers were investigated under this ``regularized transport'' framework. The work by Kurose, Yoshizawa, Amari \cite{kurose2021optimal} in this issue extended the regularizer from using Tsallis entropy to using the general deformed-entropy, a standard move adopted by researchers of deformation theory since Naudts \cite{NJ04}. The work by Tsuisui \cite{tsutsui2021optimal} in this issue investigated a generic strictly concave function as the regularizer. Explicit expressions were obtained in these works, allowing solutions to be obtained iteratively with computational efficiency. 


\subsubsection{Computational optimal transport}

One main advantage to considering the entropy-regularized problem is that it is possible to solve it much more quickly than the unregularized problem. For a more complete reference on this subject, we refer the reader to the recent text of Peyr\`e and Cuturi \cite{peyre2019computational}.

As mentioned previously, a major insight of Kantorovich was that the dual problem is a linear program, which can be solved using the simplex algorithm (or a similar such algorithm for linear programming). In practice, the convergence to the solution of an optimal transport may happen very slowly. In the worst case scenario, algorithms for solving linear programs may have exponential complexity (although their average complexity is generally polynomial).

To avoid slow convergence, one often does not try to solve such a problem exactly, but instead tries to find an approximate solution which is close to being optimal. Using an auction algorithm, one can find an arbitrarily good approximation in polynomial time. More precisely, for $X$ and $Y$ finite sets (say, both of size $N$), one can find an algorithm which terminates in $O(N^3 \log \|C\|_\infty )$ steps, where $\|C\|_\infty$ is the $L^\infty$ norm of the cost function (\cite{bertsekas1998network}, page 264).

Using the entropy-regularized optimal transport, it is possible to do even better. In particular, Cuturi also showed that solution \eqref{opt_sol} to the entropy-regulated problem $W_\epsilon(p,q)$ can be obtained by the Sinkhorn algorithm. In this algorithm, the $\alpha$ and $\beta$ vectors need not be evaluated; only the vectors $u$ and $v$, which are iteratively calculated. Sinkhorn's algorithm (with entropy of size $\epsilon$) can be solved with time complexity 
$O\left(N^{2} \log^4 (N) \epsilon^{-3}\right)$, which saves a factor of (almost) $N$ compared to the auction algorithm \cite{peyre2019computational,altschuler2017near}.

To explain why regularizing speeds up the convergence, it is helpful to understand the behavior of a non-regularized optimal transport problem. For instance, in many settings one will use the discretized optimal transport as an approximation for a continuous optimal transport problem. If the associated continuous problem satisfies the hypotheses of Theorem \ref{Deterministictransport}, the solution to the continuous problem will be a Monge solution, which the Kantorovich solution will approximate in a weak sense. As a result, the solution will concentrate on a sparse region of the product space $X \times Y$ (in this context, a set of size $O(N)$). This phenomena makes solving such problems difficult, since most points in $X \times Y$ will not be in the support of the optimal coupling. 

On the other hand, the solution to the entropy-regularized problem will be more diffuse, which makes it more amenable to numeric computation. Along with faster convergence for the computational algorithms, this phenomena also leads to faster convergence in a statistical sense (page 59 \cite{peyre2019computational}).


\subsection{Extending Legendre duality to $\lambda$-duality}

Legendre duality plays a crucial role in information geometry as it provides the canonical divergence in a dually flat space ((for exponential or mixture families of probability density functions, for instance)). In optimal transport, a more generalized notion of $c$-duality appears, based on the notion of $c$-conjugation 
\eqref{c-conjugatedef}, which generalizes the Legendre conjugation, in which $c(x,y)$ is simply to be taken as $ - \langle x, y \rangle$. 

In a series of papers, Pal and Wong \cite{PW16,PW18,W19,PW18b} applied optimal transport to mathematical finance by invoking a special cost function $c(x,y) = \log (1- \langle x, y \rangle)$. A new divergence function, called a {\it logarithmic divergence} \cite{PW16}, was proposed for a suitably generalized convex function $\psi$ that deforms the Bregman divergence. Wong \cite{W18} then extended this to a one-parameter family of cost functions
$\tilde{c}(x,y) = \kappa_\lambda(\langle x, y \rangle)$, where $\lambda$ is a real number\footnote{Originally \cite{W18} used $\alpha$ and treated $\alpha>0$ and $\alpha <0$ cases separately. This is unified by \cite{WZ21} based on the notion of  $\lambda$-exponential convexity.}  indexing the single-variable function $\kappa_\lambda(\cdot)$:
\[ \kappa_\lambda(t)=\frac{1}{\lambda} \log(1 + \lambda t) , \]
and $\log a$ is defined to be $-\infty$ if $a \leq 0$. The work of \cite{W18} showed that this generates a class of deformed exponential families for which the underlying space has \emph{constant curvature}. A one-parameter family of divergence, called $\lambda$-divergence was found which includes the logarithmic divergence ($\lambda = -1$) and Bregman divergence ($\lambda = 0$) as special cases \cite{W18,WZ21}:
\begin{eqnarray} \label{eqn:lambda.logdivergence}
D_{\lambda, \psi}(x , x') &=& \psi(x) - \psi(x') - \kappa_\lambda \left( \langle \mathsf{D} \psi(x') , x - x' \rangle \right) \\
&=& \psi(x) - \psi(x') - \frac{1}{\lambda} \log \left(1 + \lambda \langle \mathsf{D} \psi(x'), x - x' \rangle \right) \nonumber \\
&=& 
\psi (x) - \psi (x^\prime) - \kappa_\lambda \left( \langle x , u^\prime \rangle \right) + \kappa_\lambda \left( \langle x^\prime , u^\prime \rangle \right) \label{eqn:lambda.logdivergence2}
\end{eqnarray}
where $\mathsf{D} \psi$ is the usual gradient operation on $\psi$,  $\mathsf{D}^\lambda \psi$ is the deformed gradient (appropriate for the $\lambda$-deformed duality):
$$ 
\mathsf{D}^\lambda \psi(x) = \frac{\mathsf{D} \psi(x)}{1 - \lambda \langle x , \mathsf{D} \psi(x) \rangle }   ,
$$
and 
$$
u' = \mathsf{D}^{\lambda} \psi(x'), u = \mathsf{D}^\lambda \psi (x) \,\, \Longleftrightarrow \,\, x' = \mathsf{D}^{\lambda} \psi^{\lambda}(u'), x = \mathsf{D}^\lambda \psi^\lambda (u) .
$$
Here and below, $\psi^{\lambda}$ is the $\lambda$-conjugate function of $\psi$ defined, as a special case of  \eqref{c-conjugatedef}, by
$$
\psi^\lambda (u) = \sup_x \left( \frac{1}{\lambda} \log (1+\lambda \langle x, u \rangle) - \psi(x) \right) .
$$
When the $\lambda$-divergence is expressed in the conjugate variables $u, u^\prime$, a generalized version of reference--representation biduality results:
$$
D_{\lambda, \psi^{\lambda}}(u' , u)= D_{\lambda, \psi}(x , x') .
$$
This result generalizes the reference--representation biduality characteristic of Bregman divergence and the canonical divergence for dually flat spaces (see \cite{ZJ05,ZJ13}). 

In the above theory, the $\kappa_\lambda$ function associated to $c$-duality, which originates as a cost function, is viewed as a deformation to the standard Legendre duality. This $\lambda$-deformation theory, fully developed in Zhang and Wong \cite{zhang2022lambda}, enables $\lambda$-deformed Legendre duality and the standard Legendre duality to be mutually transformable to each other based on a reparameterization of one of the $\lambda$-conjugate variables $x,u$. This fact underlies the finding of \cite{WZ21}, where a $\lambda$-deformed exponential family has two apparent expressions, i.e., one with subtractive normalization and the other with divisive normalization. The former corresponds to the $q$-deformed exponential family associated to the Tsallis entropy and divergence, whereas the latter corresponds to the deformed exponential families studied by \cite{W18} associated to the R\'eyni entropy and divergence. The $\lambda$-deformation framework, inspired by the $c$-duality in optimal transport specializing the function form of $\kappa_\lambda$, appears to provide a {\it canonical} extension of the dually flat geometry to an underlying manifold with {\it constant curvature} \cite{W18,WZ21,zhang2021handbook,zhang2022lambda}. 

\subsection{Relating Otto calculus to information geometry}
\label{Works in Otto calculus/IG}

The $2$-Wasserstein geometry is special in that it admits a formal Riemannian structure (in addition to being a distance space), so one can apply Otto calculus (to the upstairs space) in connection with statistical structure (modeling the downstairs manifold both in the finite and infinite dimensional settings). For the terminology "downstairs manifold" and ``upstairs space,'' see Section \ref{otto_cal}.

\subsubsection{Upstairs and downstairs probability measures}
\label{Upstairs and downstairs section}
\begin{center}
\begin{figure}
    \centering
    \includegraphics[width=.7\linewidth]{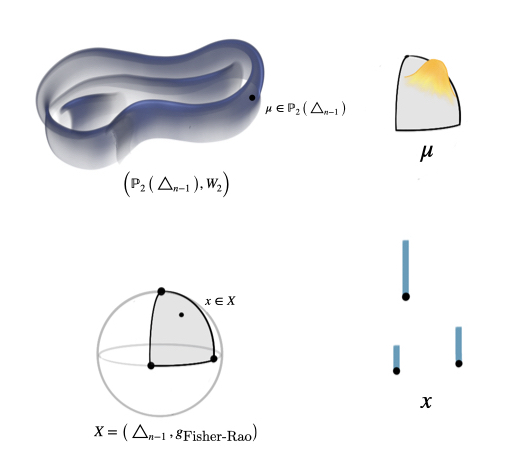}
    \caption{Upstairs and downstairs probability measures}
    \label{fig:Upstairs and downstairs}
\end{figure}    
\end{center}
One natural setting for exploring the relationship between information geometry and Otto calculus is to consider the case when the downstairs manifold $X$ is itself a statistical manifold instead of simply a Riemannian manifold (as was considered in Otto's work). A statistical manifold is a Riemannian manifold (with Fisher-Rao metric) along with a pair of conjugate connections; the statistical structure can be induced by
a proper divergence function, such as KL divergence. 

For instance, if we consider discrete probability measures on a finite set $S= \{1,\ldots,n\}$, the space of all such probability measures is the probability simplex $\triangle_{n-1}$. Here, each point on $\triangle_{n-1}$ is a discrete measure (a mixture of point measures). Taking the entirety of discrete measures supported on $S$ to be the space $X$, the Fisher-Rao metric (on the probability simplex) gives the Riemannian metric of the space of discrete measure (probability mixture). 
The statistical manifold of $\triangle_{n-1}$ is isometric to the positive orthant of an $(n-1)$-sphere, and is a well-known space in information geometry. 

Next, we take $\triangle_{n-1}$ to be the downstairs manifold, and consider probability measures $\rho$ supported on $\triangle_{n-1}$ (in Bayesian language, this is the hyper-prior on discrete probability measures). As such, we denote the space of all such probability measures by $\mathbb{P}_2(\triangle_{n-1})$; this forms the upstairs space. In this example, the Fisher-Rao metric (of the downstairs manifold) has bounded radius, so all probability measures (on the upstairs space) will have finite second-moment. Figure \ref{fig:Upstairs and downstairs} depicts both the upstairs space $\mathbb{P}_2(\triangle_{n-1})$ and the downstairs manifold $\triangle_{n-1}$. 

In linking Otto calculus with information geometry, one must be mindful that an element in the (infinite-dimensional) upstairs space $\mathbb{P}_2(\triangle_{n-1})$ should not be confused with an element in the (finite-dimensional) downstairs space
$\triangle_{n-1}$, even though both represent probability measures. However, one can embed $\triangle_{n-1}$ isometrically in $\mathbb{P}_2(\triangle_{n-1})$ via the the map $x \mapsto \delta_x$. 


For a more general statistical manifold $X$, it is necessary to restrict one's attention to $\mathbb{P}_2(X)$, so that the Wasserstein metric \eqref{Wasserstein domain} is well-defined, but the essential idea is the same. We can then study optimal transport on $\mathbb{P}_2(X)$ where the underlying cost is the distance-squared induced by the Fisher-Rao metric. It is worth noting that in the original work of Pal and Wong \cite{PW18}, a very similar framework was adopted, except that the cost for the transport was taken to be the underlying divergence, rather than the squared-distance.

In a novel approach to statistical structure of a manifold, 
Khesin et al.\ \cite{khesin2013geometry} studied the space of densities on a Riemannian manifold $\mathcal{M}=X$ in the context of its diffeomorphism group $Diff(X)$. More precisely, they treat the space of densities on $X$ as the right cosets arising from the quotient $Diff(X)/Diff_0(X)$, where $Diff_0(X)$ is the subgroup of diffeomorphisms of $\mathcal{M}$ that preserve its Riemannian volume. Notably, this work derived the Fisher-Rao metric and the Amari-Chentsov tensor on $\mathcal{M}$ in this setting. As such, this framework demonstrates a remarkably deep link to information geometry through intrinsic geometric objects such as infinite-dimensional Lie groups. It should be noted that this density manifold approach also links to optimal transport (through studying hydrodynamics of incompressible fluids, i.e., Euler equation as geodesic equation on $Diff_0(X)$), and provides metrics which interpolate between the Fisher-Rao, Sobolev $H^{1}$, and Wasserstein metrics by parametric variation). See also the work of Bauer, Bruveris, and Michor \cite{bauer2016uniqueness} who further demonstrated the uniqueness of the Fisher-Rao metric and the Amari-Chentsov tensor on this density manifold $\mathcal{M}$. 

There have been other work studying metrics which interpolate between the Wasserstein metric and Fisher-Rao geometry in the continuous setting \cite{chizat2018interpolating,kondratyev2016new,liero2018opt}. This line of work has applications in functional analysis and partial differential equations, and is an active area of research. This work is closely related to the entropy-regularized optimal transport in the continuous setting, since the Fisher-Rao metric is induced by the second variation of the relative entropy, whereas the Wasserstein metric is induced by optimal transport.

\subsubsection{Transport information geometry}

In contrast to the previous way of relating 2-Wasserstein geometry and information geometry, where the 2-Wasserstein geometry of the upstairs space relates to a statistical structure downstairs (following Otto calculus), another approach is to consider the pullback of the formal Riemannian metric to the {\it downstairs} statistical manifold $X$. This idea first appeared in McCann's work on displacement interpolation \cite{mccann1997convexity}, where he observed that the statistical manifold of normal distributions is geodesically convex in $\mathbb{P}_2(\mathbb{R}^n)$. This geometry was studied in detail by Takatsu \cite{takatsu2011wasserstein}. In related work, she also studied $q$-Gaussian measures \cite{takatsu2012wasserstein} and $\phi$-exponentials \cite{takatsu2013behaviors} using this approach.

In the past five years, a related line of work was advanced by Li and his colleagues, \cite{li2018natural,li2020ricci,chen2020optimal,li2021transport2, li2021transport,lee2021tropical}, under the name \emph{transport information geometry}. Specifically, this work compares the Fisher-Rao metric with its dual connections and the 2-Wasserstein metric ``side-by-side.'' In the case where the underlying space is a finite set with a distance $d$, the statistical manifold $X=\triangle_{n-1}$ is the probability simplex, and the distance $d$ can be used to define the 2-Wasserstein metric on $X$. This is a case of particular interest because it is analytically trackable, so we mention it now.


For instance, Li and Mont\'{u}far \cite{li2018natural} used transport information geometry to construct a natural gradient flow, which they called the Wasserstein gradient flow, using the Otto's prescription of the Riemannian metric. In this issue, Li \cite{li2021transport2} studies the Bregman divergence and KL divergence for transport information geometry in the case where the underlying space is finite. 
The readers may wish to compare and contrast the above approach with the entropy-regularized transport framework (see Section \ref{Entropy-regularized transport}). 

When this approach is applied to the space of measures on a finite graph $G$ with a distance function $d$, interesting results were obtained when comparing the associated Wasserstein distance (induced by the graph metric) with the Fisher-Rao metric, see \cite{li2021transport} of this issue.

Due to the differing topologies induced by the Wasserstein geometry and Fisher-Rao metric (see Section \ref{Differing topologies}), comparing the Wasserstein and Fisher-Rao geometries in the infinite-dimensional setting presents several challenges. To circumvent these issues, Chen and Li \cite{chen2020optimal} consider statistical manifolds $\mathcal{P}_\theta$ as finite-dimensional submanifolds of $\mathbb{P}_2(X)$. By pulling back Otto's Riemannian structure to $\mathcal{P}_\theta$, this induces a Riemannian metric associated to optimal transport. It is worth noting that except for the case of normal distributions and a few other special cases, the manifold $\mathcal{P}_\theta$ will not be totally geodesic in $\mathbb{P}_2(X)$, so the minimizing paths in $\mathcal{P}_\theta$ will not be the same as those in $\mathbb{P}_2(X)$.

\subsection{Other papers included in this Special Issue}

Previously, 
 Malago, Montrucchio and Pistone  \cite{malago2018wasserstein} studied $2$-Wasserstein geometry of multivariate normal distributions (i.e., Gaussian measures) and carried out explicit calculations of Riemannian geodesics and Levi-Civita covariant derivatives over the space of symmetric matrices, which parametrizes the Gaussian measure as an exponential model. In this issue, Montrucchio and Pistone \cite{montrucchio2021kantorovich} investigated the $1$-Wasserstein distance for finite metric space and representation of $l_1$-embeddable metrics. They performed in-depth analysis of weighted trees and tree-like spaces. This work contains novel results about optimal transport on trees, e.g., deriving the closed-form expression for the $1$-Wasserstein transport between measures on a tree in term of cut norms, and showing that the cut norm is equivalent to the Arens-Eells norm in this setting.

 Mallasto, Gerolin, and Minh \cite{mallasto2021entropy} in this issue investigated entropy-regularized 2-Wasserstein distance between Gaussian measures. Through involved matrix calculations, they provided a closed-form expression of this distance. It is well-known that the $2$-Wasserstein metric between Gaussians is the Bures metric comparing quantum states in quantum information geometry. By incorporating a regularization parameter, the above authors are able to compare it with the Sinkorn divergence as its limiting case. Their work is among several other concurrent groups that independently provided calculations about this important case of Gaussian measure transport (e.g., \cite{takatsu2011wasserstein}). 

Finally, Sei \cite{sei2021coordinate} in this Special Issue studied multi-dimensional probability distributions and showed that under certain regularity conditions, there exists a unique coordinate-wise transformation such that the transformed distribution satisfies a Stein-type identity. These transformations naturally arise in copula theory modeling joint distributions as having uniform marginals. The main results  of that paper include a variational characterization of the possible target measures, a proof of uniqueness of such maps, and a sufficient condition (namely, copositivity) for their existence. The proofs used some arguments from optimal transport theory; it was based on an energy minimization problem over a totally geodesic subset of the Wasserstein space. Sei's result is considered as an alternative to Sklar's theorem regarding copulas, and is also interpreted as a
generalization of a diagonal scaling theorem.

\section{K\"ahler versus para-K\"ahler Geometries for Optimal Transport}
\label{Geometrizing Optimal Transport section}
A fruitful interaction between information geometry and optimal transport stems from attempts to prescribe a ``canonical'' geometry to optimal transport problems. There are multiple ways to geometrize optimal transport, and we will mention two such theories which have close links with information geometry.

\subsection{The para-K\"ahler geometry of optimal transport}

While studying the regularity theory of the Monge problem, Kim and McCann \cite{kim2010continuity} found a natural geometric formulation of optimal transport in terms of para-K\"ahler geometry. More precisely, given two smooth $n$-dimensional manifolds $X$ and $Y$ with Borel measures $\mu$ and $\nu$ and a cost function $c:X \times Y \to \mathbb{R}$ which satisfies several regularity and non-degeneracy conditions,
 they considered the space $X \times Y$ as a pseudo-Riemannian manifold with the pseudo-metric
\begin{equation}\label{pseudoR}
h:=\frac{1}{2}\left(\begin{array}{cc}
0 & -c_{i,j} \\
-c_{i,j} & 0
\end{array}\right).
\end{equation}
Here $c_{i,j} = \partial_{x^i} \partial_{y^j} c(x,y)$. This induces the product space with a pseudo-metric of signature $(n,n)$ (i.e., $n$ positive eigenvalues and $n$ negative ones).
 It turns out \cite{cortes2004special,villani2008stability} that this is also an example of a para-K\"ahler manifold, 
 which is a pseudo-Riemannian metric with an associated para-complex structure $K$ satisfying
  \[h(K\eta, K \xi) = -h (\eta, \xi). \]
 Here the para-complex structure $K$ induces the map $K (\eta,\xi) = (\eta, -\xi)$ for $\eta \in T_x X$ and $\xi \in T_y Y$, and $K^2$ equals the identity map. Para-K\"ahler metrics of the above form were discussed by the second named author in the context of information geometry \cite{zhang2016information}.
 
 
 
 To study an optimal transport problem on this space, Kim, McCann and Warren \cite{kim2010pseudo} consider the conformal deformation 
 \begin{equation}
h_{c}^{\rho, \bar{\rho}}=\left(\frac{\rho(x) \bar{\rho}(\bar{x}))}{|\operatorname{det}(c_{i,j})|}\right)^{\frac{1}{n}}\left(\begin{array}{cc}
0 & -c_{i,j} \\
-c_{i,j} & 0
\end{array}\right)
\end{equation}
 which is another pseudo-metric on $X \times Y$. In this new geometry, the solution to the Monge problem is an $n$-dimensional stable maximal surface. Furthermore, they developed a theory of calibrations for pseudo-Riemannian metrics and showed that the solution to the optimal transport problem is calibrated by the function
 \begin{equation}
\Phi(x, \bar{x})=\frac{\rho(x) d x+\bar{\rho}(\bar{x}) d \bar{x}}{2}.
\end{equation}

Wong and Yang \cite{wong2021pseudo} considered the framework of a $c$-divergence (based on \cite{PW18} and \cite{W18}) and considered the graph of the optimal transport as a (possibly disconnected) statistical manifold equipped with the c-divergence which is induced by the cost function and the pair of Kantorovich potentials. In this work, they show that the dual connections $\nabla$ and $\nabla^\star$ are the projections of the Levi-Civita connection of the pseudo-Riemannian metric \eqref{pseudoR} on the product space. This gives an embedding of the dualistic information geometry into pseudo-Riemmanian geometry, and demonstrates a deep connection between optimal transport and information geometry.

Another striking feature of the para-K\"ahler framework of optimal transport is that the curvature of the product space encodes the MTW condition. In the Kim-McCann geometry, the MTW tensor is the curvature of certain light-like planes. This interpretation gives intrinsic differential geometric structure to the regularity problem and immediately explains many of the properties of $\mathfrak{S}$, such as why it transforms tensorially under change of coordinates. We refer the reader to (\cite{kim2010continuity}, Section 5) for some interesting remarks on how various aspects of optimal transport can be understood by studying the geometry of this space.
 
\subsection{K\"ahler geometry of optimal transport}

For cost functions which are induced by a convex potential $\Psi$ (such as $c(x,y)=\Psi(x-y)$ or a $\mathcal{D}_\Psi^{(\alpha)}$-divergence \cite{zhang2004divergence}), there is a separate formulation of the regularity theory in terms of K\"ahlerian information geometry \cite{khan2020kahler}. In particular, given such a cost function, one can use the convex potential to define a K\"ahler metric whose curvature encodes the MTW tensor.



For this, we consider a convex potential $\Psi:\Omega \to \mathbb{R}$, where
 $\Omega$ is assumed to contain the relevant Minkowski linear combination of the sets $X$ and $Y$ so that the cost function is well defined.

When one treats the domain $\Omega$ as a Hessian manifold (with a global chart), its tangent bundle $T\Omega$ admits a canonical K\"ahler structure which has a natural $\mathbb{R}^n$ symmetry induced by translations in the fiber directions.\footnote{It is a general fact that the tangent bundle of any Hessian manifold admits a K\"ahler metric (see, e.g. \cite{dombrowski1962geometry,satoh2007almost}).}

In order to relate the MTW tensor to the curvature of the tangent bundle $T\Omega$, it is necessary to lift tangent vectors from $\Omega$ to corresponding vectors in $T \Omega$. For this, we use the \emph{polarized lift}.

\begin{definition}
Given a tangent vector $v \in T_{x_0} \Omega$ defined as
\begin{equation}
v =\sum_i  v^i \frac{\partial}{\partial {x^i}}, \end{equation}
we define its polarized lift to be the $(1,0)$ vector 
\begin{equation} 
v^{pol.} := \sum_i v^i \frac{\partial}{\partial {z^i}} \in T_{(x_0,y)}^{(1,0)}T \Omega.
\end{equation}

\end{definition}

Here, $\{ z^i\}_{i=1}^n$ are the holomorphic coordinates $x^i + \sqrt{-1} y^i$. Since the metric is translation symmetric in the fibers, we are free to pick the fiber coordinate $y$ at will.  These vectors are polarized in the sense that their components are real (in the $z$-coordinates).
By considering polarized $(1,0)$ vectors, we can make sense of the notion of positive or negative anti-bisectional curvature.

\begin{definition}[Anti-bisectional curvature]
 Given two $(1,0)$ vectors $V$ and $W$, the anti-bisectional curvature is defined to be
\begin{equation}
     \mathfrak{A}(V,W)= R(V,\overline W, V, \overline W).
\end{equation}
\end{definition}

At first, it does not seem meaningful to discuss positive and negative anti-bisectional curvature; this is because  $\mathfrak{A}(\sqrt{-1}V,W) = - \mathfrak{A}(V,W)$.
 However, by restricting our attention to polarized vectors, it is meaningful to talk about this quantity as having a sign, and it turns out to be proportional to the MTW tensor.

\begin{theorem}[\cite{khan2020kahler}] \label{Kahlergeometryofoptimaltransport}
For a cost which is induced by a convex potential $\Psi: \Omega \to \mathbb{R}$, the MTW tensor $\mathfrak{S}(\xi,\eta)$ is proportional\footnote{The proportionality constant depends on how the convex potential is used to induce the cost function.} to the anti-bisectional curvature of an associated K\"ahler-Sasaki metric on the tangent bundle $T \Omega$.
\begin{equation}
  \mathfrak{S}(\xi,\eta) \propto   R \left( \xi^{pol.}, \overline{ (\eta^\sharp)^{pol.}},\xi^{pol.}, \overline{ (\eta^\sharp)^{pol.}} \right) .
\end{equation}  
\end{theorem}
 Furthermore, the conjugate connection (of the associated Hessian manifold $\Omega$) encodes the notion of relative $c$-convexity.
 
 \begin{proposition} \label{Relativecconvexity}
For a cost of the form $c(x,y) = \Psi(x-y)$, a set $Y$ is $c$-convex relative to $X$ if and only if, for all $x \in X$, the set $x-Y$ is geodesically convex with respect to the dual connection $\nabla^*.$ 
\end{proposition}
 
At present, this geometry only applies to cost functions which are induced by a convex potential $\Psi$. However, this includes some important examples. In particular, using this geometry 
the present authors \cite{khan2020kahler} were able to establish a regularity result for relative pseudo-arbitrages, which solved a question posed by Pal and Wong in their work \cite{PW18b} relating information geometry, optimal transport, and mathematical finance.

\subsubsection{Some open problems related to this geometry}
One open problem is to extend this work to a non-K\"ahler theory for costs defined on a Lie group $G$. More precisely, on any Lie group $G$ we can construct a curvature-free connection $\nabla_G$ whose torsion is induced by the associated Lie algebra. The connection can then be used to define an almost complex structure on the tangent bundle $T G$, which is non-integrable whenever $G$ is non-Abelian. Then, given a potential function $\Psi: G \to \mathbb{R}$, it may be possible to define a structure analogous to a metric on $G$ (although care must be taken since $\nabla_G^2 \Psi$ will generally not be symmetric, so does not define a Riemannian metric). If this can be done, the hope is that it can then be used to define a Sasaki-type structure on the tangent bundle $TG$ whose curvature encodes the MTW tensor of the cost $c(x,y) =\Psi(xy^{-1})$. From a high level perspective, the K\"ahler formulation of optimal transport is precisely this construction for the \emph{Abelian} Lie group $\mathbb{R}^n$ (and where $\Psi$ may be defined on an open subset of $\mathbb{R}^n$ rather than the entire space). 

\begin{question}
Is there a non-K\"ahler version of Theorem \ref{Kahlergeometryofoptimaltransport} for costs of the form \begin{equation}\label{Lie group psi cost}
    c(x,y) =\Psi(xy^{-1}) 
\end{equation}where $x$ and $y$ are elements of a Lie group?
\end{question}

As a particular application of the previous question, one might hope to find examples of cost functions which satisfy the MTW(0) condition. In particular, the squared distance function for a left-invariant metric on $SO(3)$ is precisely a function of the form \eqref{Lie group psi cost}, which raises the following question.

\begin{question}
Is it possible to classify left-invariant metrics on $SO(3)$ whose squared distance costs have non-negative MTW tensor?
\end{question}

Another natural question about the complex geometry of optimal transport is how to interpret the geometric meaning of the solution of an optimal transport problem. It seems that this is intrinsically linked with displacement interpolation, but it is not clear how to make this precise.

Using displacement interpolation, it is possible to understand the solution to optimal transport in terms of a section $s$ of $T X$. To do so, we consider a flow which follows the $c$-exponential map and whose initial condition is $s$. In general, this flow will satisfy a particular semi-spray equation induced by the cost function. In the setting of 2-Wasserstein geometry, the flow lines will be geodesics in the usual sense. For costs induced by a convex potential, the flow lines will be geodesics with respect to the \emph{dual} connection. Along this flow, we can solve the continuity equation with the initial measure $\mu = \mu(0)$. The section $s$ encodes a viable transport plan if the induced measure at time $t=1$ is our desired target measure $\nu$. Furthermore, this transport plan will be optimal if the induced cost from this flow is minimal. 

This gives a way to encode the solutions to optimal transport in terms of sections to the tangent bundle, but it is not clear what geometric properties such a section will need to satisfy in order to induce the optimal transport between $\mu$ and $\nu$. Furthermore, it is not clear what role the curvature of $T \Omega$ plays (which may be a different space from $T X$).

\begin{question}
When does a section of $T X$ induce an optimal flow (in the sense of displacement interpolation)? What is the relationship between these sections and the curvature of the tangent bundle $T \Omega$?
\end{question}

Recall that in the para-complex geometry, the MTW is the curvature of certain light-like planes, whereas in the complex setting, the MTW tensor is the orthogonal anti-bisectional curvature. The link between the K\"ahler geometry approach and the para-K\"ahler geometry approach is not well understood. 

\begin{question}
How can one relate the pseudo-Riemannian geometry to the complex geometry of optimal transport? 
\end{question}

Before moving on, it is worth noting that the K\"ahler manifolds which appear in this theory also have applications outside optimal transport. In particular, they exhibit a form of mirror symmetry induced by the information-geometric duality of the Hessian metric \cite{zhang2020statistical}. As shown in \cite{khan2021hall}, this duality appears to relate to the theory of automorphic forms encountered in analytic number theory.

\subsection{K\"ahler-Ricci flow and other smoothing flows}
Although the complex geometric perspective is less general than the pseudo-Riemannian framework, there are several advantages associated to using a Hermitian, rather than an indefinite, pseudo-metric (para-Hermitian metric). In particular, it is possible to use techniques from complex geometry, such as K\"ahler-Ricci flow, to deform the cost functions and study optimal transport. 

Ricci flow \cite{hamilton1982three} is a geometric flow which deforms a Riemannian metric 
by its Ricci curvature:
\begin{equation}\label{Ricciflow}
    \frac{\partial g}{\partial t}  = - 2 \, \textrm{Ric}(g).
\end{equation}
When the underlying manifold is K\"ahler, this flow preserves the complex structure and the K\"ahlerian nature of the metric \cite{cao1985deformation}, so is known as the K\"ahler-Ricci flow. Recently, the first named author and Zheng investigated the behavior of the MTW tensor under  K\"ahler-Ricci flow \cite{khan2020k}.

\begin{theorem}[\cite{khan2020k}] \label{2DKRflowpreservesMTW0}
Suppose that $ \Omega \subset \mathbb{R}^n $ is a convex domain and $\Psi :\Omega \to \mathbb{R}$ is a strongly convex function so that the associated K\"ahler manifold $(T \Omega, \omega_0)$  is complete and has bounded curvature.
\begin{enumerate}
    \item For $n = 2$, \emph{non-negative orthogonal anti-bisectional curvature} (i.e., the MTW(0) condition) is preserved by K\"ahler-Ricci flow.
    \item For all $n$, \emph{non-positive anti-bisectional curvature} is preserved under K\"ahler-Ricci flow.
\end{enumerate}
\end{theorem}

This result allows us to prove regularity for the two-dimensional Monge problem for cost functions which are induced by a convex potential and are much rougher than $C^4$, which is the regularity needed so that Equation \eqref{MTW0} is meaningful (see also, \cite{guillen2015local,loeper2021weak}). For instance, it is possible to prove a H\"older estimate on the transport for costs which are $W^{2,p}$ for $p>2$ and which satisfy the MTW condition in the sense that the positive-time potentials induce K\"ahler metrics with non-negative anti-bisectional curvature. 

\subsubsection{Smoothing flows in optimal transport}

This line of work suggests that smoothing flows can play a role in optimal transport, especially in cases with low regularity. Parabolic flows for the transport potential $\psi$ mentioned in Section \ref{sec_regularity} had previously been considered \cite{kitagawa2012parabolic} but in this context K\"ahler-Ricci flow is used to deform the \textit{cost function}.
 
  Currently, Theorem \ref{2DKRflowpreservesMTW0} is limited to costs which are induced by a convex potential because it is not possible to define a Ricci flow naively for general cost functions. The reason for this is that the Ricci flow is not weakly parabolic in the pseudo-Riemannian setting. More concretely, the natural analogue of the K\"ahler-Ricci flow equation in the Kim-McCann geometry is the following:
 \begin{equation} \label{ParaKahlerRicciflow}
     \frac{\partial}{\partial t} c(x,y) = \log \det[ c_{i,j }(x,y) ].
 \end{equation}
Solutions to this equation are not unique. For instance, given one solution  $c(x,y,t)$, there are infinitely many other solutions of the form \[ c(x,y,t) + f(x,t)+g(y,t)\]
 where $f(x,t)$ and $g(y,t)$ are arbitrary functions of $(x,t)$ and $(y,t)$, respectively. In other words, this flow cannot be smoothing, and it is not even clear when solutions exist at all. It would be of interest to find a smoothing flow for more general costs, which we leave as an open problem.
 
 \begin{question}
  For more general cost functions, is it possible to define a smoothing flow?
 \end{question}

 When the cost function is the squared distance on a Riemannian manifold, one natural idea for a smoothing flow is to evolve the underlying metric by Ricci flow (see, e.g., \cite{mccann2010ricci,topping2009L}). This is distinct from  \eqref{ParaKahlerRicciflow}, and will be smoothing away from the cut-locus. This raises the following question.

\begin{question}
For the squared-distance cost on a Riemannian manifold, does the Ricci flow preserve non-negativity of the MTW tensor?
\end{question}

We suspect that the answer to this question is positive, at least in dimensions two and three. The reason to believe this conjecture is that for the squared distance cost, the MTW tensor restricted to the diagonal (i.e., $x=y$) is two-thirds the sectional curvature. As such, the MTW(0) condition is a non-local strengthening of non-negative sectional curvature \cite{loeper2009regularity}. In two and three dimensions, the Ricci flow preserves positivity of the sectional curvature and Hamilton showed that if a metric starts with positive sectional curvature, it converges to a metric of constant curvature with the scalar curvature satisfying a differential Harnack estimate. As such, two- and three-dimensional manifolds with non-negative MTW tensor converge quickly to metrics of constant curvature. Furthermore, Figalli, Rifford and Villani showed that in all dimensions, small deformations (in the $C^4$-sense) of the round sphere satisfy the MTW condition \cite{figalli2012nearly}.

Together, these results show that if the initial metric satisfies the MTW condition, there will be a neighborhood of the diagonal for which the MTW tensor is positive for all positive time of the flow. Furthermore, after a definite amount of time, the MTW condition will hold everywhere. In higher dimensions, Ricci flow need not preserve positivity of the sectional curvature, so it might be necessary to make the additional assumption that the curvature operator is positive (and to formulate the correct version of positivity for the MTW operator). In spirit, these results can be established using Hamilton's tensor maximum principle (at least away from the cut locus), but the computations needed are extremely involved. 


\section{Conclusion}
In this paper, we have surveyed several lines of fruitful interaction between optimal transport and information geometry: entropy-regularized (also called entropy-relaxed) transport, $c$-duality as generalizing Legendre duality, relating Wasserstein metric and Otto calculus to the statistical structure,  para-K\"ahler versus K\"ahler geometry, etc. Our survey is admittedly incomplete: we have focused largely on the work published in the Springer journal \emph{Information Geometry}, including the current Special Issue, with an intended audience of information geometers. For this reason, we have omitted a large amount of work on optimal transport/Wasserstein geometry that does not have an apparent information geometry connection. 

We conclude that the budding intersection between optimal transport and information geometry is a highly promising area of research. Given that both Wasserstein geometry and information geometry have, in their own rights, seen wide applications in statistics, machine learning, computer vision, etc., it is all the more remarkable that there are some deep, internal links of these two mathematical sub-disciplines. Future research will continue to illuminate the relation between these two geometric frameworks about probability measures (and probability density functions).  

\bigskip 

{\bf Acknowledgement} The authors would like to thank Leonard Wong for his helpful comments on the manuscript. G.K. is partially supported by a Simons Collaboration Grant 849022 (``K\"ahler-Ricci flow and optimal transport"), while J.Z. is partially supported by United States Air Force Office for Scientific Research, grant number AFOSR-FA9550-19-1-0213.

\bibliography{references.bib}
\bibliographystyle{alpha}

\end{document}